\documentclass[10pt]{article}
\usepackage{amssymb,theorem}
\makeatletter
\@addtoreset{equation}{section}
\makeatother
\renewcommand{\ge}{\mathop{\geqslant}}
\renewcommand{\le}{\mathop{\leqslant}}
\newcommand{\E}{\mathbb{E}}
\renewcommand{\P}{\mathbb{P}}

\theoremheaderfont{\sffamily\bfseries\upshape}
\theoremstyle{break}{
\newtheorem{theorem}{Theorem}[section]

\newtheorem{proposition}[theorem]{Proposition}
}
\newcommand {\qed}%
{%
    {}\hfill
    {}\hfill
    {$\square $}%
    \vspace {0.3cm}%
    \pagebreak [2]%
    \par
}%
\newenvironment{proof}[1]{%
        \vspace{0.3cm}%
        \pagebreak [2]%
        \par%
        \noindent {\bf  Proof ~ \textbf{#1}\ }}{\qed}%
\newenvironment{remark}{%
        \vspace{0.3cm} \pagebreak [2]%
        \par%
        \refstepcounter{theorem}%
        \noindent%
        {\bfseries\sffamily Remark~\thetheorem}\\}{\qed}%
\newcommand{\eps}{\varepsilon}

\title{Sojourn Times of Brownian Sheet}
\author{D. Khoshnevisan\thanks{Research supported in part by
    grants from NSF and NATO}\\University of Utah
    \and R. Pemantle\thanks{Research supported in part by
    NSF grant 98-03249}\\Ohio State University
}
\date{September 26, 2000}
\begin{document}
\maketitle
\renewcommand{\abstractname}{}
\begin{abstract}
\noindent\rule{\linewidth}{.5mm}
\noindent\rule{\linewidth}{.5mm}\\
\begin{center}
\sffamily\slshape
This paper is dedicated to Professor Endr\'e Cs\'aki\\on the occasion
of his 65th birthday.
\end{center}

\end{abstract}

\section{Introduction}
Let $B$ denote the standard Brownian sheet. 
That is, $B$ is a centered Gaussian process indexed
by $\mathbb{R}^2_+$ with continuous trajectories and
covariance structure
\[
   \mathbb{E}\big\{ B_s B_t\big\} =\min\{ s_1,t_1\} 
  	\times \min\{ s_2,t_2\},\qquad
  	s=(s_1,s_2) ,\ t=(t_1,t_2)\in\mathbb{R}^2_+.
\]
In a canonical way, one can think of $B$ as ``two-parameter
Brownian motion''.

In this article, we address the following question:
``\emph{Given a measurable function 
$\upsilon:\mathbb{R}\to\mathbb{R}_+$,
what can be said about the distribution of 
$\int_{[0,1]^2} \upsilon(B_s)\ ds$?}''
The one-parameter variant of this question is both
easy-to-state and well understood. Indeed, if 
$b$ designates standard Brownian motion, the Laplace transform
of $\int_0^1 \upsilon(b_s+x)\ ds$ often solves a Dirichlet 
eigenvalue problem (in $x$), as prescribed by the Feynman--Kac
formula; cf. Revuz and Yor \cite{RevuzYor}, for example.
While analogues of Feynman-Kac for $B$ are not yet known to hold,
the following highlights some of the unusual behavior of 
$\int_{[0,1]^2}\upsilon(B_s)\ ds$ in case 
$\upsilon=\mathbf{1}_{[0,\infty)}$ 
and, anecdotally, implies that finding
explicit formul\ae\ may present a challenging task.

\begin{theorem}
   \label{thm:sojourn}
   There exists a $c_0\in(0,1)$, such that
   for all $0<\varepsilon<\frac{1}{8}$,
   \[
	\exp\Big\{- \frac{1}{c_0}\log^2(1/\varepsilon) \Big\}
	\le
	\mathbb{P}\big\{ 
		{\textstyle\int}_{[0,1]^2}\mathbf{1}_{\{ B_s>0\}}
		\ ds <\varepsilon \big\} \le
		\exp\Big\{ -c_0{
		\log^2(1/\varepsilon)}\Big\}.
   \]
\end{theorem}

\begin{remark}
	\label{rem:OneParArcSin}%
	By the arcsine law, the one-parameter version of the above has
	the following simple form: given a linear Brownian motion $b$,
	\[
	   \lim_{\eps\to 0^+} \eps^{-1/2} \P\big\{ 
	   {\textstyle\int}_0^1\mathbf{1}_{\{ b_s>0\}}\ ds
	   <\eps \big\} =\frac {2}{\pi};
	\]
	see \cite[Theorem 2.7, Ch. 6]{RevuzYor}.
\end{remark}

\begin{remark}
	R. Pyke (personal communication) has asked whether
	$\int_{[0,1]^2}\mathbf{1}{\{ B_s>0\}}\ ds$ has an
	arcsine-type law; see \cite[Section 4.3.2]{Pyke}
	for a variant of this question in discrete time.
	According to Theorem \ref{thm:sojourn},
	as $\eps\to 0$, the cumulative distribution function of
	$\int_{[0,1]^2}\mathbf{1}_{\{ B_s>0\}}\ ds$ goes to zero
	faster than any power of $\eps$. In particular, 
	the distribution of time (in $[0,1]^2$) spent positive
	does not have any simple extension of the arsine law.
\end{remark}

\begin{theorem}
   \label{thm:smallball}
   Let $\upsilon (x) := \mathbf{1}_{[-1,1]}(x)$, 
   or $\upsilon (x) := \mathbf{1}_{(-\infty,1)}(x)$. Then,
   there exists a $c_1\in(0,1)$, such that for all $\eps\in(0,\frac{1}{8})$,
   \[
	\exp\Big\{ -{\log^3(1/\eps)\over c_1 \eps}\Big\}\le
	\mathbb{P}\big\{ {\textstyle\int}_{[0,1]^2} 
	\upsilon (B_s)\ ds <\eps\big\}
	\le
	\exp\Big\{ -c_1{\log(1/\eps)\over \eps}\Big\}.
   \]
\end{theorem}

For a refinement, see Theorem \ref{thm:refsmallball} below.

\begin{remark}
	\label{rem:OneDimFK}%
	The one-parameter version of Theorem \ref{thm:smallball} is quite simple. 
	For example, let $\Gamma=\int_0^1\mathbf{1}_{[-1,1]}(b_s)\ ds$,
	where $b$ is linear Brownian motion. In principle, one can compute the Laplace
	transform of $\Gamma$ by means of Kac's formula
	and invert it to calculate its distribution function.
	However, direct arguments suffice to show that the two-parameter Theorem \ref{thm:smallball}
	is more subtle than its one-parameter counterpart:
	\begin{equation}
	\label{eq:OneDimFK}
		-\infty <
		\liminf_{\eps\to 0^+} \eps \ln \P \{ \Gamma<\eps\} \le
		\limsup_{\eps\to 0^+} \eps \ln \P \{ \Gamma<\eps\} < 0,
	\end{equation}
	where $\ln$ denotes the natural logarithm function.
	We will verify this later on in the Appendix.
\end{remark}

\begin{remark}
	The arguments used to demonstrate Theorem \ref{thm:smallball} 
	can be used to also estimate the distribution function
	of additive functionals
	of form, e.g., $\int_{[0,1]^2}\upsilon(B_s)\ ds$, as long as 
	$\alpha \mathbf{1}_{[-r,r]}\le \upsilon \le \beta \mathbf{1}_{[-R,R]}$, where
	$0<r\le R$ and $0<\alpha\le\beta$. Other formulations are also possible. 
	For instance, when 
	$\alpha \mathbf{1}_{(-\infty,r)}\le \upsilon\le 
	\beta\mathbf{1}_{(-\infty,R)}.$
\end{remark}

\section{Proof of Theorems \ref{thm:sojourn} and \ref{thm:smallball}}

Our proof of Theorem \ref{thm:sojourn} rests on a
lemma that is close in spirit to a Feynman--Kac formula 
of the theory of one-parameter Markov processes.

\begin{proposition}
   \label{prop:FK}
   There exists a finite and positive constant $c_2$, such that
   for all measurable $D\subset\mathbb{R}$ and
   all $0<\eta,\varepsilon<\frac{1}{8}$.
   \[
     	\mathbb{P}\big\{ {\textstyle\int}_{[0,1]^2} \mathbf{1}_{\{
     		B_s \not \in D\}}\ ds <\varepsilon\big\}
      	\le \mathbb{P}\big\{ \forall s\in[0,1]^2:~
     		B_s  \in D_{\varepsilon^{\frac{1}{4}-2\eta}}\big\}+
    		\exp\{ -c_2\eps^{-\eta}\},
   \]
   where $D_\delta$ denotes the $\delta$-enlargement of $D$ for
   any $\delta>0$. That is, 
   \[ 
       	D_\delta := \big\{ x\in\mathbb{R}:~ \mathrm{dist}
    			( x; D)\le\delta\big\},
   \]
   where `$\mathrm{dist}$' denotes Hausdorff distance. 
\end{proposition}

\begin{proof}{}
For all $t\in[0,1]^2$, let $|t| := \max\{ t_1, t_2\}$. Then,
it is clear that for any $\eps,\delta>0$,
whenever there exists some $s_0\in[0,1]^2$ for which
$B_{s_0}\not\in D_\delta$, either
\begin{enumerate}
	\item $\sup_{|t-s|\le\eps^{1/2}} |B_t-B_s|>\delta$, where
    		the supremum is taken over all such choices of
    		$s$ and $t$ in $[0,1]^2$; or
	\item for all $t\in[0,1]^2$ with $|t-s_0|\le\eps^{1/2}$,
    	$B_t \in D$, in which case, we can certainly deduce that
     $\int_{[0,1]^2} \mathbf{1}_{D^\complement}(B_t)\ dt>\eps.$
\end{enumerate}
Thus,
\begin{eqnarray*}
	\mathbb{P}\big\{
	\exists s_0\in[0,1]^2:~ B_{s_0} \not \in D_\delta \big\} &\le &
	\mathbb{P}\big\{ \textstyle \sup_{|t-s|\le\eps^{1/2}} 
	|B_t-B_s|>\delta\big\} + \\
	&& \qquad + \mathbb{P}\big\{ {\textstyle\int}_{[0,1]^2}
		\mathbf{1}_{D^\complement}(B_t)\ dt >\eps\big\}.
\end{eqnarray*}
By the general theory of Gaussian processes, there exists a 
universal positive and finite constant $c_2$ such that
\begin{equation}
	\label{eq:maximal}
	\mathbb{P}\big\{ \sup_{|t-s|\le\eps^{1/2}} 
	|B_t-B_s|>\delta\big\}  \le \exp\big\{
	-c_2\delta^2\eps^{-1/2}\big\}.
\end{equation}
Although it is well known, we include a brief derivation of this 
inequality for completeness.
Indeed, we recall C. Borell's inequality from 
Adler \cite[Theorem 2.1]{Adler}: 
if $\{ g_t;~ t\in T\}$ is a centered Gaussian process
such that $\|g\|_T = \E\{\sup_{t\in T}|g_t|\}<\infty$ and whenever
$T$ is totally bounded in the metric $d(s,t)=\sqrt{\E\{ (g_t-g_s)^2\}}$
($s,t\in T$),
\[
   \P\{ \sup_{t\in T} |g_t| \ge \lambda + \|g\|_T \} 
   \le 2 \exp\Big\{ -{\lambda^2\over 2\sigma^2_T}\Big\},
\]
where $\sigma_T^2 =\sup_{t\in T}\E\{ g_t^2\}$.
Eq. (\ref{eq:maximal}) follows from this by letting
$T=\{ (s,t)\in(0,1)^2\times(0,1)^2:~|s-t|\le\eps^{1/2}\}$,
$g_{t,s} = B_t-B_s$ and by making a few lines of standard calculations.
Having derived (\ref{eq:maximal}), 
we can let $\delta := \eps^{\frac{1}{4}-\frac{\eta}{2}}$ to obtain
the proposition.
\end{proof}

\noindent\textbf{Proof of Theorem \ref{thm:sojourn}}\
Let $D=(-\infty,0)$ and use Proposition \ref{prop:FK} to see that 
\[
   \mathbb{P}\big\{
   {\textstyle\int}_{[0,1]^N}
   \mathbf{1}_{\{ B_s>0\}}<\varepsilon\big\}
   \le \mathbb{P}\big\{ \sup_{s\in [0,1]^2}B_s\le\varepsilon^{
   {1\over 4}-2\eta}\big\} + \exp\{-c_2\varepsilon^{-\eta}\}.
\]
Thus, the upper bound of Theorem \ref{thm:sojourn} follows from 
Li and Shao \cite{LiShao}, which states that
\[
   \limsup_{\eps\to 0^+} {1\over\log^2(1/\eps)}\
   \log\P\big\{ \sup_{s\in[0,1]^2} B_s \le \eps\big\}<-\infty.
\]
(An earlier, less refined version, of this estimate can be found  in  Cs\'aki et al. \cite{CKS}.)
To prove the lower bound, we note that
\begin{eqnarray*}
   \lefteqn{\mathbb{P}\big\{ {\textstyle\int}_{[0,1]^2}\mathbf{1}_{\{ B_s>0\}}\
	ds < 2\varepsilon-\eps^2 \big\}}\\
   && \ge  \mathbb{P}\big\{ 
	\sup_{s\in[\varepsilon,1]^2} B_s < 0 \big\}\\
   && =  \mathbb{P}\big\{ \forall (u,v)\in [ 
		0,\ln ( \textstyle{1\over\varepsilon} )]^2:~
		e^{(u+v)/2} \; B(e^{-u},e^{-v}) <0\big\},
\end{eqnarray*}
and observe that the stochastic process
$(u,v)\mapsto B( e^{-u},e^{-v})/ e^{-(u+v)/2}$ is the
2-parameter Ornstein--Uhlenbeck sheet. 
All that we need
to know about the latter process
is that it is a stationary, positively correlated Gaussian process whose
law is supported on the space of continuous functions on $[0,1]^2$.
We define $c_3>0$ via the equation
\[
   e^{-c_3} := \mathbb{P}\Big\{ \forall (u,v)\in[0,1]^2:~
	{B(e^{-u},e^{-v})\over e^{-(u+v)/2}}<0\Big\}.
\]
By the support theorem, $0<c_3<\infty$; this is a consequence
of the Cameron-Martin theorem on Gauss space; cf.
Janson \cite[Theorem 14.1]{Janson}. Moreover,
by stationarity and by Slepian's inequality (cf. \cite[Corollary 2.4]{Adler}),
\begin{eqnarray*}
   \lefteqn{\mathbb{P}\big\{ {\textstyle\int}_{[0,1]^2}
	\mathbf{1}_{\{ B_s<0\}}\
	ds < \eps \big\} }\\
   &&\ge
	\prod_{0\le i,j\le \ln(1/\varepsilon)+1}
   	\mathbb{P}\Big\{ \forall (u,v)\in[i,i+1]\times[j,j+1]:~
		{B(e^{-u},e^{-v})\over e^{-(u+v)/2}}<0\Big\}\\
   && = \exp\Big\{ -c_3 \ln^2(e^2/\varepsilon)\Big\}.
\end{eqnarray*}
This proves the theorem.\hfill$\square$\\

Next, we prove Theorem \ref{thm:smallball}.\\

\noindent\textbf{Proof of Theorem \ref{thm:smallball}}\
Let $\mathcal{D}_\eps$ denote the collection of all points
$(s,t)\in[0,1]^2$, such that $st\le \eps$.
Note that 
\begin{enumerate}
	\item Lebesgue's measure of $\mathcal{D}_\eps$ is at least
		$\eps\ln(1/\eps)$; and
	\item if $\sup_{s\in \mathcal{D}_\eps} |B_s|\le 1$, then
		$\int_{[0,1]^2}\mathbf{1}_{(-1,1)}(B_s)\ ds >\eps\ln(1/\eps)$.
\end{enumerate}
Thus,
\[
	\mathbb{P}\Big\{ \int_{[0,1]^2}\mathbf{1}_{(-1,1)}(B_s)\ ds 
	< \eps\ln(1/\eps)\Big\}
	\le \mathbb{P}\Big\{ \sup_{s\in\mathcal{D}_\eps} |B_s| >1\Big\}.
\]
A basic feature of the set $\mathcal{D}_\eps$ is that whenever 
$s\in\mathcal{D}_\eps$, then $\mathbb{E}\{ B_s^2\}\le\eps$. 
Since
$\E\{ \sup_{s\in\mathcal{D}_\eps}|B_s|\} \le \E\{ \sup_{s\in[0,1]^2}|B_s|\}<\infty$,
we can apply Borell's inequality to deduce the existence of a finite, positive
constant $c_4<1$, such that for all $\eps>0$,
\(
	\mathbb{P}\{ \sup_{s\in\mathcal{D}_\eps} |B_s| >1/c_4\}
	\le \exp\{ -{c_4/\eps}\}.
\)
We apply Brownian scaling and possibly adjust $c_4$ to conclude that
\[
	\mathbb{P}\Big\{ \sup_{s\in\mathcal{D}_\eps} |B_s| >1\Big\}
	\le e^{-c_4/\eps}.
\]
Consequently, we can find a positive, finite constant $c_5$, such that
for all $\eps\in(0,\frac{1}{8})$,
\begin{equation}
	\P\{\Gamma< \eps\} \le \exp\Big\{ -c_5{\ln(1/\eps)\over\eps}\Big\}.
	\label{eq:lower-smallball}
\end{equation}
This implies the upper bound in the conclusion of Theorem \ref{thm:smallball}.
For the lower bound, we note that for all 
$\eps\in(0,\frac{1}{8})$, Lebesgue's measure of 
$\mathcal{D}_\eps$ is bounded above by $c_6\eps\log(1/\eps)$.
Thus,
\[
	\mathbb{P}\Big\{ \int_{[0,1]^2}\mathbf{1}_{(-\infty,1)}(B_s)\ 
	ds < c_6\eps\log(1/\eps)\Big\}
	\ge \mathbb{P}\Big\{ \inf_{s\in[0,1]^2\setminus
	\mathcal{D}_\eps} B_s >1\Big\}.
\]
On the other hand, whenever $s\in[0,1]^2\setminus\mathcal{D}_\eps$, 
$s_1s_2\ge\eps$. Thus,
\begin{eqnarray*}
	\mathbb{P}\Big\{ \int_{[0,1]^2}\mathbf{1}_{(-\infty,1)}(B_s)\ 
	ds < c_6\eps\log(1/\eps)\Big\}
	& \ge & \mathbb{P}\Big\{ \inf_{s\in[0,1]^2\setminus
		\mathcal{D}_\eps} {B_s\over
		\sqrt{s_1s_2}} >{1\over\sqrt{\eps}}\Big\}\\
	& = & \mathbb{P}\Big\{ \inf_{\scriptstyle u,v\ge 0:\atop\scriptstyle
			u+v\le\ln(1/\eps)} O_{u,v} >\eps^{-1/2}\Big\},
\end{eqnarray*}
where $O_{u,v} := B(e^{-u},e^{-v})/ e^{-(u+v)/2}$ is an 
Ornstein--Uhlenbeck sheet.
Consequently,
\[
	\mathbb{P}\Big\{ \int_{[0,1]^2}\mathbf{1}_{(-\infty,1)}(B_s)\ 
	ds < c_6\eps\log(1/\eps)\Big\}
	\ge \mathbb{P}\Big\{ \inf_{0 \le u,v\le \ln(1/\eps)} 
	O_{u,v} >\eps^{-1/2}\Big\},
\]
By appealing to Slepian's inequality and to the stationarity of $O$, 
we can deduce that
\begin{eqnarray}
	\lefteqn{\mathbb{P}\Big\{ \int_{[0,1]^2}\mathbf{1}_{(-\infty,1)}
	(B_s)\ ds < c_3\eps\log(1/\eps)\Big\}} \qquad\qquad\nonumber \\
	&& \ge  \prod_{0\le i,j\le \ln(1/\eps)}  \mathbb{P}\Big\{
	\inf_{i \le u\le i+1}\:\inf_{j\le v\le j+1} 
			O_{u,v} >\eps^{-1/2}\Big\} \nonumber \\
	&& = \bigg[ \mathbb{P}\Big\{ \inf_{0 \le u,v\le 1} O_{u,v}
		>\eps^{-1/2}\Big\} \bigg]^{ \ln^2(e/\eps)}.
	\label{eq:lower}
\end{eqnarray}
On the other hand, recalling the construction of $O$, we have
\begin{eqnarray*}
	\lefteqn{\mathbb{P}\Big\{ \inf_{0\le u,v\le 1} O_{u,v} > 
	\eps^{-1/2}\Big\} }\\
	&&\ge  \mathbb{P}\Big\{ \inf_{1\le s,t\le e} B_{s,t} \ge
		e\ \eps^{-1/2}\Big\} \\
	&& \ge  \mathbb{P}\Big\{ B_{1,1} \ge 2e\ \eps^{-1/2}\ ,\
		\sup_{1\le s_1,s_2\le e}\big| B_s-B_{1,1}\big|\le 
		e\ \eps^{-1/2}\Big\}\\
	&& = \mathbb{P}\Big\{ B_{1,1} \ge 2e\ \eps^{-1/2}\Big\}\cdot
		\mathbb{P}\Big\{ \sup_{1\le s_1,s_2\le e}\big| 
		B_s-B_{1,1}\big|\le e\ \eps^{-1/2}\Big\}\\
	&& \ge c_7 \mathbb{P}\Big\{ B_{1,1} \ge 2e\ \eps^{-1/2}\Big\},
\end{eqnarray*}
for some absolute constant $c_7$ that is chosen independently of all 
$\eps\in(0,\frac{1}{8})$.
Therefore, by picking $c_8$ large enough, we can insure that for all 
$\eps\in(0,\frac{1}{8})$,
\[
	\mathbb{P}\Big\{ \inf_{0\le u,v\le 1} O_{u,v} >\eps^{-1/2}\Big\}
	\ge \exp\big\{ -c_8\eps^{-1}\big\}.
\]
Plugging this in to Eq. (\ref{eq:lower}), we obtain
\begin{equation}
	\label{eq:upper-smallball}
	\mathbb{P}\Big\{ \int_{[0,1]^2}\mathbf{1}_{(-\infty,1)}(B_s)\ ds
			< c_6\eps\log(1/\eps)\Big\} \ge 
		\exp\Big\{ -c_8{\ln^2(1/\eps)\over 4 \eps}\Big\}.
\end{equation}
The lower bound of Theorem \ref{thm:smallball}
follows from replacing $\eps$ by $\eps/\ln(1/\eps)$.\hfill$\square$\\

The methods of this proof go through with few changes to 
derive the following extension of Theorem \ref{thm:smallball}.

\begin{theorem}
   \label{thm:refsmallball}
   Suppose $\varphi:\mathbb{R}_+\to\mathbb{R}_+$ is a measurable function
   such that (a) as $r\downarrow 0$, $\varphi(r)\uparrow\infty$; and
   (b) there exists a finite constant $\gamma>0$, such that 
   for all $r\in(0,\frac{1}{2})$, $\varphi(2r)\ge \gamma\varphi(r)$.
   Define $J_\varphi = \int_{[0,1]^2} \mathbf{1}_{\{ |B_s|\le
	\sqrt{s_1s_2}\varphi(s_1s_2)\}}\ ds.$ Then, there exist
   a finite constant $c_9>1$, such that for all 
   $\eps\in(0,\frac{1}{2})$,
   \[
	\exp\Big\{ -c_9 \varphi^2({\eps\over\log(1/\eps)}) \log^2(1/\eps)\Big\} \le
 	\mathbb{P}\big\{ J_\varphi <\eps \big\}
	\le \exp\Big\{ -\frac{1}{c_9} \varphi^2({\eps\over\log(1/\eps)})
		\Big\}.
   \]
\end{theorem}

\section*{Appendix: On Remark \ref{rem:OneDimFK}}

In this appendix, we include a brief verification of the exponential
form of the distribution function of $\Gamma$; cf. Eq. (\ref{eq:OneDimFK}). 
Given any $\lambda>\frac{1}{2}$ and for $\zeta=(2\lambda)^{-1/2}$, we have
\begin{eqnarray}
	\E\{ e^{-\lambda\Gamma} \}
	& \le & \E\Big\{ \exp \big( -\lambda{\textstyle\int}_0^\zeta \upsilon(b_s)\ ds \big) \Big\} \nonumber \\
	& \le & e^{-\lambda\zeta} + \P\{ \sup_{0\le s\le \zeta}|b_s|>1\}\\
	& \le & e^{-\lambda\zeta} + e^{-1/(2\zeta)} \nonumber \\
	& = &  2 e^{-\sqrt{\lambda/2}}.\label{eq:Laplace1}%
\end{eqnarray}
By Chebyshev's inequality,
$\P\big\{ {\textstyle\int}_0^1 \upsilon (b_s)\ ds<\eps\big\}
\le 2 \inf_{\lambda>1}e^{-\sqrt{\lambda/2} + \lambda\eps}.$
Choose $\lambda=\frac{1}{8}\eps^{-2}$ to obtain the following
for all $\eps\in(0,\frac{1}{2})$:
\begin{equation}
	\label{eq:Laplace}
	\P\big\{ \Gamma <\eps\big\}
	\le 2 e^{-1/(8\eps)}.
\end{equation}
Conversely, we can choose $\delta=(2\lambda)^{-1/2}$ and $\eta\in(0,\frac{1}{100})$
to see that
\begin{eqnarray*}
	\E\{ e^{-\lambda\Gamma}\}
	& \ge & \E\Big\{ \exp\big( -\lambda{\textstyle\int}_0^\delta \upsilon(b_s)\ ds \big) ;~
	\inf_{\delta\le s\le 1} |b_s| >1 \Big\}\\
	& \ge & e^{-\lambda\delta} \ \P\big\{ 
	|b_\delta|>1+\eta ~,~ \sup_{\delta <s<1+\delta} |b_s-b_\delta|<\eta\big\}.
\end{eqnarray*}
Thus, we can always find a positive, finite constant $c_{10}$ that only depends on
$\eta$ and such that
\[
	\E\{ e^{-\lambda\Gamma} \}
	\ge c_{10} \exp\Big\{ -\sqrt{\frac{\lambda}{2}}\ \big[ 1+(1+\eta)^2(1+\psi_\delta)\big]\Big\},
\]
where $\lim_{\delta\to 0^+}\psi_\delta=0$, uniformly in $\eta\in (0,\frac{1}{100})$.
In particular, after negotiating the constants, we obtain
\begin{equation}
	\label{eq:Laplace2}%
	\liminf_{\lambda\to\infty} \lambda^{-1/2} \ln \E\{ 
	e^{-\lambda\Gamma} \} \ge  -2^{1/2}.
\end{equation}
Thus, for any $\eps\in(0,\frac{1}{100})$,
\[
	e^{-\sqrt{2\lambda}( 1+o_1(1) )}  \le  \E\{ e^{- \lambda \Gamma}\}
	 \le  \P\big\{ \Gamma<\eps \big\} + e^{-\lambda\eps},
\]
where $o_1(1)\to 0$, as $\lambda\to\infty$, uniformly in $\eps\in(0,\frac{1}{100})$.
In particular, if we choose $\eps = (1+\eta)\sqrt{2/\lambda}$, where $\eta>0$,
we obtain
\[
	\P\big\{ \Gamma < (1+\eta)\sqrt{2/\lambda} \big\} \ge 
	e^{-\sqrt{2\lambda}( 1+o_2(1))},
\]
where $o_2(1)\to 0$, as $\lambda\to\infty$. This, Eq. (\ref{eq:Laplace})
and a few lines of calculations,
together imply Eq. (\ref{eq:OneDimFK}).\hfill$\square$\\

\begin{small}

\end{small}

\bigskip

\noindent\begin{small}
\noindent\begin{tabular}{ll}
\textsc{Davar Khoshnevisan} & \textsc{Robin Pemantle} \cr
University of Utah & Ohio State University \cr
Department of Mathematics & Department of Mathematics \cr
155 S 1400 E JWB 233 & 231 W. 18 Ave., Columbus, OH 43210 \cr
Salt Lake City, UT 84112--0090 & Columbus, OH 43210 \cr
\texttt{davar@math.utah.edu} & \texttt{pemantle@math.ohio-state.edu} \cr
\texttt{http://www.math.utah.edu/\~{}davar} &
\texttt{http://www.math.ohio-state.edu/\~{}pemantle}
\end{tabular}
\end{small}

\end{document}